# A GPU accelerated mixed-precision Finite Difference informed Random Walker (FDiRW) solver for strongly inhomogeneous diffusion problems

**Zirui Mao[1]\*, Shenyang Hu[1], Ang Li[1,2]\***

[1] Pacific Northwest National Laboratory, Richland, WA 99354
[2] University of Washington, Seattle, WA 98195

zirui.mao@pnnl.gov; ang.li@pnnl.gov

**Abstract**

In nature, many complex multi-physics coupling problems exhibit significant diffusivity inhomogeneity, where one process occurs several orders of magnitude faster than others in temporal. Simulating rapid diffusion alongside slower processes demands intensive computational resources due to the necessity for small time steps. To address these computational challenges, we have developed an efficient numerical solver named Finite Difference informed Random Walker (FDiRW). In this study, we propose a GPU-accelerated, mixed-precision configuration for the FDiRW solver to maximize efficiency through GPU multi-threaded parallel computation and lower precision computation. Numerical evaluation results reveal that the proposed GPU-accelerated mixed-precision FDiRW solver can achieve a 117X speedup over the CPU baseline, while an additional 1.75X speedup by employing lower precision GPU computation. Notably, for large model sizes, the GPU-accelerated mixed-precision FDiRW solver demonstrates strong scaling with the number of nodes used in simulation. When simulating radionuclide absorption processes by porous wasteform particles with a medium-sized model of $192 \times 192 \times 192$, this approach reduces the total computational time to 10 minutes, enabling the simulation of larger systems with strongly inhomogeneous diffusivity.

## 1. Introduction

Many natural systems exhibit strongly inhomogeneous diffusion processes, where one process occurs at a significantly faster rate than the other concurrent processes. For example, in nuclear materials, the long-term radiation decay of fuel materials coexists with the rapid evolution of radiation-induced defects across polycrystalline structures[1,2]. In γU-Mo nuclear fuels, the diffusivity of uranium (U) along the grain boundary can be 3 to 8 orders of magnitude higher than its self-diffusivity. For molybdenum (Mo), this difference ranges from 4 to 10 orders of magnitude. Similarly, in the context of radionuclide absorption by porous wasteform materials within a flowing waste stream, the diffusivity of species in solid and liquid phases can differ by 3 to 8 orders of magnitude[3,4]. These extreme variations in diffusivity present a substantial challenge in modeling such systems due to intensive computational demands.

Such massive computation demands arise from the discretization concept adopted in most mature numerical methods, such as the finite difference method (FDM)[5-7] and the finite volume method[8-10] in computational fluid dynamics (CFD) modeling. In these methods, the material or domain is discretized into a finite number of infinitesimally small divisions (elements, cells, nodes, or particles). The time marching interval, $\Delta t$, must be sufficiently small such that the physical phenomena occurring in each small division can be viewed as linear and do not have to propagate beyond the cell domain[11]. Consequently, when the diffusivity spans a wide range, the smallest diffusivity controls the long diffusion time, while the largest diffusivity determines the smallest time step $\Delta t$ required in numerical simulations. This can result in an enormous number of iterations, often reaching as many as $10^{10}$ or more, leading to prohibitively high computational burden.

To address this heavy computational burden, we have proposed a superposition-based numerical solver called the Finite Difference integrated Random Walk (FDiRW)[12,13]. This solver is designed to simulate fast diffusion process and is integrated with traditional discretization-based numerical methods to handle slower processes. In traditional discretization-based methods, long-time, long-range fast diffusion is treated as iterative time-marching of short-time, short-range diffusion. In contrast, the FDiRW solver adopts a superposition concept inspired by the Random Walk model[14-16], where the long-time, long-range process is treated as the superposition of long-time, long-range diffusion results sourcing from each point. Unlike the Random Walk model, which derives long-time, long-range diffusion results directly from analytical expression[14-16] and works only for isotropic pure diffusion, the FDiRW solver performs Finite Difference (FD) simulations to obtain the long-time, long-range diffusion results sourcing from each point. This FD simulation step is necessary because the analytical expressions adopted in Random Walk model cannot handle the anisotropic diffusion caused by the complex diffusion domain in porous structures.

The FDiRW has two obvious advantages over traditional FDM in computing fast diffusion. First, it is not established upon discretization concepts, and thus free from the limitation for small time-steps. Second, it eliminates the need for small grid sizes, because the mass distribution after a long-time fast diffusion within domain would become rather smooth, negating the need for fine mesh resolution. Consequently, the FDiRW solver can use large time steps and coarse meshes, greatly enhancing computational efficiency compared to discretization-based methods. Our previous study[12] demonstrated the FDiRW solver's potential by achieving a 1100x speedup compared to the Finite Difference solver for a 192×192×192 model. When simulating the strongly inhomogeneous diffusion processes with an integrated framework, the computational time fraction spent on solving fast diffusion was significantly reduced from 99.8% using the FD solver to 34% using the FDiRW solver.

This study investigates further efficiency enhancements of the FDiRW solver by employing GPUs for multi-thread parallel computation and exploring the emerging low-precision computation on GPUs. Specifically, the accuracy of the low-precision based FDiRW solver is examined carefully. To balance between accuracy and efficiency, we propose a GPU-accelerated mixed-precision FDiRW solver in this study. The efficiency and accuracy of the proposed GPU-accelerated mixed-precision based FDiRW solver are evaluated and compared to other benchmark configurations using CPUs and high-precision data types. Finally, the scalability of the GPU-accelerated mixed-precision based FDiRW solver to larger systems is studied by applying it to models of varying sizes. This comprehensive study aims to showcase the performance potential of the FDiRW solver in terms of accuracy and efficiency for handling computationally challenging strongly inhomogeneous diffusion problems in engineering, with a specific focus on nuclear waste stream treatment as an example.

Section 2 presents the integrated model for simulating the strongly inhomogeneous diffusion processes during nuclear waste stream treatment. Section 3 introduces the GPU-accelerated, mixed-precision FDiRW solver, and its accuracy and efficiency are investigated in Section 4. In Section 5, the proposed FDiRW solver is applied to various model sizes to evaluate its computational complexity when scaled up to larger systems.

## 2. Integrated model of nuclear waste stream treatment

During nuclear waste stream treatment, the absorption of radionuclide ions by porous wasteform particles involves a multi-physics process with strongly inhomogeneous diffusivity[17-19]. The process entails immersing micron-level porous wasteform particles into a nuclear waste solution containing radionuclide ions. These radionuclide ions diffuse rapidly in the liquid phase and are absorbed by the particle surface, then enter the solid phase, where they undergo significantly slower diffusion — typically exhibiting diffusion rates 3 to 8 orders of magnitude slower than in the liquid phase[3,4]. Accurately modeling the kinetics of radionuclide absorption by the wasteform particles is crucial for understanding the absorption mechanism and designing high-performance wasteform materials for nuclear waste stream treatment.

Considering the unaffordable computational demand caused by the strongly inhomogeneous nature of such systems, two simplifications are typically employed to alleviate the computational burden. Firstly, the volume ratio of the waste stream to the wasteform particles is assessed. The computational focus is directed towards a single wasteform particle surrounded by the corresponding volume of waste stream. This approach effectively approximates the radionuclide absorption kinetics of the entire system by studying the behavior of a single particle without compromising the underlying physics. Secondly, recognizing that the kinetics of radionuclide absorption are primarily influenced by the distribution of radionuclides in the liquid phase surrounding the wasteform particles, the computational domain is subdivided into near-field and far-field regions. This division is delineated by defining a sphere slightly larger than the wasteform particle, as represented by the dashed circle in Figure 1. Suppose the considered particle has a radius of $r_P$. We define the liquid phase within the sphere with a radius of $r_p + 5\Delta h$ as the near-field liquid phase, with $\Delta h$ being the grid size. In computational analyses, only the fast diffusion in the near-field domain is considered, while the radionuclide concentration in the far-field is assumed to be uniformly distributed. This approach significantly reduces computational burden by focusing on the region most crucial to the absorption process while still considering the overall system behavior.

Overall, the computational process involves four main components: (1) the fast diffusion of radionuclide ions in the near-field liquid; (2) the absorption of radionuclide ions at the particle surface; (3) the slow diffusion of radionuclide ions within the solid particle, and (4) the calculation of the radionuclide concentration remaining in the far-field derived from concentration conservation based on the total concentration within the solid and near-field liquid, as shown in Figure 1.

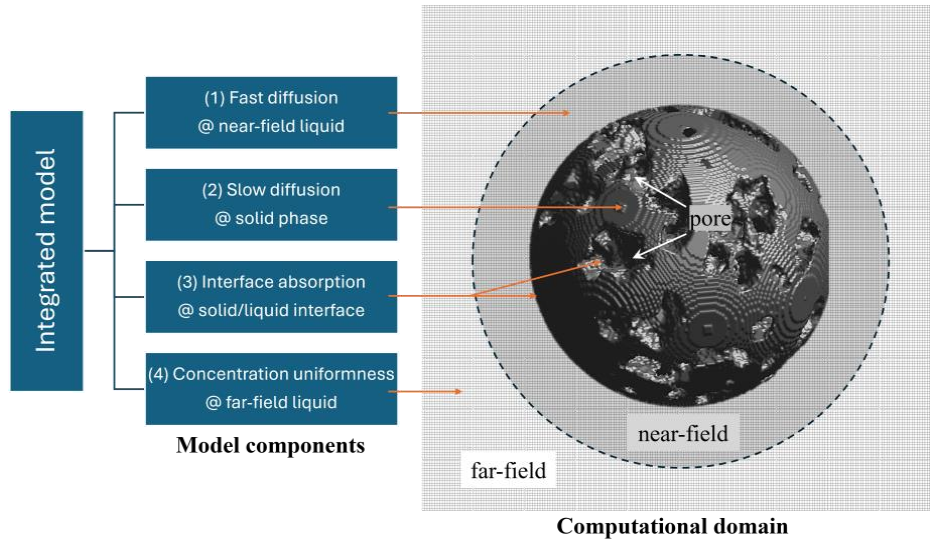

Figure 1 Computational domains surrounding the porous wasteform particle. The porous wastform particle is generated numerically with phase field modeling by referring to the microstructure information measured in experiment[20].

The multi-physics process is described by the diffusion-reaction equation, which has been validated in our previous work[12,21,22]:

$$\frac{\partial c_m(\mathbf{r}, t)}{\partial t} = \nabla \cdot \left( D_m \nabla \left( \frac{\mu_m}{RT} \right) \right) + \dot{R}_m, \qquad m = \text{S in } solid \text{ and } m = \text{L in } liquid\ phase, \tag{1}$$

Here, $\mathbf{r} = (x, y, z)$ is the spatial coordinate and $t$ is time, $c_m$ is the instantaneous radionuclide ion concentration in the phase $m$, and $D_m$ is the diffusivity of radionuclide ions in phase $m$. $R$ and $T$ are the gas constant and temperature in K, respectively. $\mu_m$ is the chemical potential governing the diffusion behavior in the liquid phase and absorption capacity of the solid phase. $\dot{R}_m$ is the absorption rate on the solid/liquid interface.

*Diffusion in solid and liquid phases in near-field*. The first term on the right-hand side of Eq. (1) describes ion diffusion, driven by the gradient of chemical potential $\mu_m$. The chemical potential drives the ion concentration $c_m$ in phase $m$ towards its equilibrium value $c_m^{eq}$ through defining chemical potential as

$$\mu_S = \mu_S(\mathbf{r}, t) = A_S\big(c_S(\mathbf{r}, t) - c_S^{eq}\big) \tag{2}$$

$$\mu_L = \mu_L(\mathbf{r}, t) = A_L\big(c_L(\mathbf{r}, t) - c_L^{eq}\big) \tag{3}$$

in the solid and liquid phases, respectively. Here, $A_S$ and $A_L$ are the free energy coefficients for the solid phase and liquid phases, respectively. $c_S^{eq}$ is the equilibrium concentration in the solid phase, representing the maximum concentration of radionuclides absorbed by the solid particle as measured in experiments, while $c_L^{eq}$ is the equilibrium concentration in the liquid phase at the final equilibrium state [13,21]. Initially, the concentration of radionuclide in the solid phase is lower than $c_S^{eq}$, resulting in a negative chemical potential. In contrast, the liquid phase has a higher initial concentration than $c_L^{eq}$, thus yielding a positive chemical potential. The higher chemical potential in the liquid phase compared to the solid phase drives the radionuclides to penetrate from the liquid phase into the solid phase, which simulates the ions absorption process from the liquid phase to the solid phase.

*Interface absorption*. The absorption rate of radionuclides at interface is governed by the interface reaction rate $\dot{R}_m$, *i.e.*, the second term in the right-hand side of Eq. (1). This rate follows a Pseudo-Second-Order reaction mechanism[20,23] represented by Eqs. (4)-(6).

$$\dot{R}_S = kf_L f_S, \ \dot{R}_L = -\dot{R}_S \tag{4}$$

$$f_L = \begin{cases} 0, & as\ c_L(t) \leq c_L^{eq} \\ \frac{c_L(t) - c_L^{eq}}{c_L^{eq}}, & as\ c_L(t) > c_L^{eq} \end{cases} \tag{5}$$

$$f_S = \frac{c_S^{eq} - c_S(t)}{c_S^{eq}}. \tag{6}$$

Basically, the interface reaction occurs only when the ion concentration in the liquid exceeds the equilibrium value, *i.e.*, $c_L(t) > c_L^{eq}$, and when the concentration in the solid is below its equilibrium value, *i.e.*, $c_S(t) < c_S^{eq}$. With those conditions, the final absorption rate is controlled by a constant $k$ in Eq. (4).

*Far-field concentration*. When modeling such diffusion-reaction systems, it is essential to configure proper boundary conditions, specifically for the domain outside the near-field domain, known as the far-field domain. Since the far-field concentration is assumed instantly uniform, the radionuclide concentration in the far-field can be updated according to the concentration conservation constraint at every time step,

$$c_L^{far-field}(t) = \frac{\Sigma c_{S+L}(t_0) - \Sigma c_L^{near-field}(t) - \Sigma c_S(t)}{V_L^{far-field}} \tag{7}$$

Here, $\Sigma c_{S+L}(t_0)$ is the total concentration of radionuclides in the whole domain including solid phase and liquid phase initially. $\Sigma c_L^{near-field}(t)$ and $\Sigma c_S(t)$ are the instantaneous total concentration in the near-field liquid and solid, respectively. $V_L^{far-field}$ is the volume of far-field liquid. The averaged uniform concentration $c_L^{far-field}$, calcualted by Eq. (7), conserves the total concentration in the entire domain.

All the necessary parameter values and material properties are summarized in Table 1. The initial and equilibrium concentrations in solid and liquid phases are configured to ensure concentration conservation between the initial and equilibrium states, based on the volume fractions of liquid and solid and experimental data[24]. The effectiveness of the integrated model associated with the listed parameters in simulating nuclear waste stream treatment have been verified in our previous works [12,13,20].

Table 1 Parameter values and material properties employed in simulation for the case of particle size $r_P = 50\Delta h$.

| $\Delta h$ | $10 \times 10^{-9}\ [m]$ | $\frac{A_L}{RT}$ | $2 \times 10^3$ |
|---|---|---|---|

| $\Delta t$ | 500 $[\mu s]$ | $c_S^{eq}$ | 1.0 |
|---|---|---|---|
| $\Delta t_{fd}$ | 0.5 $[\mu s]$ | $c_L^{eq}$ | $1.0 \times 10^{-5}$ |
| $D_S$ | $1.0 \times 10^{-17} \ [m^2/s]$ | $c_S^0$ | $1.0 \times 10^{-6}$ |
| $D_L$ | $1.0 \times 10^{-14} \ [m^2/s]$ | $c_L^0$ | $2.12 \times 10^{-3}$ |
| $\frac{A_S}{RT}$ | $2 \times 10^3$ | $k$ | 0.05 $[1/s]$ |
| $V_{Solid}$ | $3.67 \times 10^{-13}$ [mL] | $V_{Liquid}$ | $1.83 \times 10^{-10}$ [mL] |
| $V_L^{far-field}$ | $1.80 \times 10^{-10}$ [mL] | $\Sigma c_{S+L}(t_0)$ | 388,716.10 |

## 3. FDiRW solver for fast diffusion in near-field liquid

### 3.1 Coarse mesh-based FDiRW algorithm for fast diffusion

For an isolated diffusion process, the mass proportion moving between locations remains constant within a fixed long time period. The basic idea of the FDiRW solver is to restore these preconditioned proportions as $\mathbf{p}$ and treat the fast diffusion problems across the domain as the superposition of these proportions, weighted by the instantaneous concentration at the corresponding locations:

$$\mathbf{c}(t + \Delta t) = \mathbf{p}(\Delta t)\mathbf{c}(t_0) + \mathbf{p}_{BC}(\Delta t)c_{far-field}(t) \qquad (8)$$

with

$$\mathbf{c} = \begin{bmatrix} c_1 \\ c_2 \\ \vdots \\ c_{N_L} \end{bmatrix}, \mathbf{p} = \begin{bmatrix} p_{11} & p_{12} & \cdots & p_{1N_L} \\ p_{21} & p_{22} & \cdots & p_{2N_L} \\ \vdots & \vdots & \ddots & \vdots \\ p_{N_L1} & p_{N_L2} & \cdots & p_{N_LN_L} \end{bmatrix}, \mathbf{P}_{BC} = \begin{bmatrix} p_1^{BC} \\ p_2^{BC} \\ \vdots \\ p_{N_L}^{BC} \end{bmatrix} \qquad (9)$$

where $\mathbf{c}$ is the concentration array of $N_L$ liquid nodes, $\mathbf{p}_{BC}$ is the proportions list corresponding to the effect of boundary conditions, and $c_{far-field}$ is the concentration at boundaries.

In the coefficient matrix $\mathbf{p}$, $p_{ij}$ is the mass proportion moving from node $j$ to node $i$. The elements of $j$-th column in $\mathbf{p}$ can be obtained directly by solving the governing diffusion equation with an initial single point source at node $i$ using the explicit Finite Difference Method. Considering that the concentration distribution after long-time fast diffusion is rather smooth in spatial, a coarse mesh is sufficient to represent the mass proportions list, offering more promising efficiency. Accordingly, we proposed the coarse mesh-based FDiRW[12,25]. To distinguish from concentration array $\mathbf{c}$ and coefficient matrix $\mathbf{p}$ corresponding to full nodes, the coarse mesh resultant ones are denoted as $\mathbf{C}$ and $\mathbf{P}$, respectively, with reduced dimensions. Thus, the coarse mesh-based FDiRW solver is expressed as

$$\mathbf{C}(t + \Delta t) = \mathbf{P}(\Delta t)\mathbf{C}(t) + \mathbf{P}_{BC}(\Delta t)c_{far-field}(t_0) \qquad (10)$$

A mesh-coarsening algorithm and variable-mapping algorithm between the fine and coarse meshes were developed and verified in our previous work[12] to generate the exclusive coarse mesh and connect the array $\mathbf{C}$ on the coarse mesh to the array $\mathbf{c}$ on the fine mesh. With that, the total number of nodes involved in the FDiRW computation can be greatly reduced from $N_L$ to $N$, where $N$ is the number of nodes in the coarse mesh. For a balance between accuracy and efficiency, in this study, we used a coarse mesh with $N = \frac{1}{125} N_L$. In the 3-D case, this means the average size of the coarse mesh is five times that of the fine mesh.

Therefore, the fast diffusion in the liquid phase can be solved straightforwardly with Eq. (10) once the coefficient matrices are obtained using the given method. Note that before performing the FDiRW calculation, the concentration matrix $\mathbf{C}$ in the coarse mesh must be calculated based on the concentration matrix $\mathbf{c}$ in the fine mesh, referred as the 'mapping' step. After the FDiRW calculation, the concentration matrix $\mathbf{c}$ in fine mesh must be computed based on the $\mathbf{C}$ matrix obtained from FDiRW, referred as the reverse 're-mapping' step.

In detail, the 'mapping' step converting $\mathbf{c}$ matrix corresponding to the fine mesh into $\mathbf{C}$ matrix corresponding the coarse mesh, is expressed as:

$$C(t, \mathbf{r}_I) = \frac{\sum_{i=1}^{N_I} c(t, \mathbf{r}_i)(\Delta h)^d}{\sum_{i=1}^{N_I} (\Delta h)^d} \tag{11}$$

whereas the 'remapping' step reversely converting $\mathbf{C}$ matrix into $\mathbf{c}$ matrix is:

$$c(t + \Delta t, \mathbf{r}_i) = C(t + \Delta t, \mathbf{r}_I) \tag{12}$$

That is, in the 'mapping' step, the averaged concentration of the liquid nodes is assigned to the corresponding representative node. After solving the fast diffusion in the liquid phase, the liquid nodes belonging to the same group are assumed to have a uniform concentration. Although the 'mapping' and 'remapping' steps represented by Eqs. (11) and (12) are simple, they explicitly conserve the concentration, making them suitable for use in such applications.

For more details about the FDiRW approach and mesh-coarsening algorithm, please refer to our previous work[12]. The three-step coarse mesh-based FDiRW solver can be represented as:

'Mapping' step: $C(t, \boldsymbol{r}_I) = \frac{\sum_{i=1}^{N_I} c(t, \mathbf{r}_i)(\Delta h)^d}{\sum_{i=1}^{N_I} (\Delta h)^d} = \frac{\sum_{i=1}^{N_I} c(t, \mathbf{r}_i)}{N_I}$ for each representative node $I$ $\tag{13}$

FDiRW calculation: $C(t + \Delta t, \boldsymbol{r}_I) = \sum_{i=1}^{N} \{P_{1,i} * C_i\} + c_{far-field} * P_1^{BC}$ for each representative node $I$ $\tag{14}$

'Re-mapping' step: $c(t + \Delta t, \boldsymbol{r}_i) = C(t + \Delta t, \boldsymbol{r}_I)$ for each liquid node $i$ $\tag{15}$

### 3.2 GPU acceleration of FDiRW algorithm

In the three-step FDiRW solver represented by Eqs. (13)-(15), computations at each node are performed independently, making FDiRW naturally suitable to parallel computation on multi-thread GPUs. In this study, we propose a GPU-accelerated FDiRW solver, as illustrated in Figure 2, to achieve optimal efficiency. Note that the three-step FDiRW algorithm is inherently sequential, with each step dependent on the preceding one. However, the method can still benefit substantially from parallel processing within each step.

In the 'mapping' and 'remapping' steps, the calculations for all liquid nodes are evenly distributed across 256 GPU threads to enable parallel computation. We avoid distributing computational loads based on representative liquid nodes because each representative node in the coarse mesh may correspond to a varying number of liquid nodes in the fine mesh. This variability can lead to unbalanced workloads, causing some GPU threads to be overburdened and slowing down the overall process. By assigning tasks based on fine-mesh liquid nodes, we ensure the same workloads for each thread, optimizing computational efficiency. For the FDiRW calculation step, the matrices $\mathbf{P}$ and $\mathbf{P}_{BC}$ are distributed evenly across 256 GPU threads. Each GPU thread handles a submatrix of $\mathbf{P}$ and multiplies it concurrently by the entire matrix $\mathbf{C}$. This approach divides the matrix-vector multiplication into 256 smaller, concurrent submatrix-vector multiplications, thus enhancing efficiency. The number of threads is set to 256 to achieve maximum efficiency, as increasing the number of threads beyond this does not yield any further improvements.

We wrote a subfunction to perform matrix-vector multiplication in CUDA. Meanwhile, various CUDA libraries, such as cuBLAS, offer optimized matrix-vector multiplication routines that can enhance performance. In this study, we also examined the efficiency gains achievable by employing the optimized cuBLAS package in GPU parallel computation.

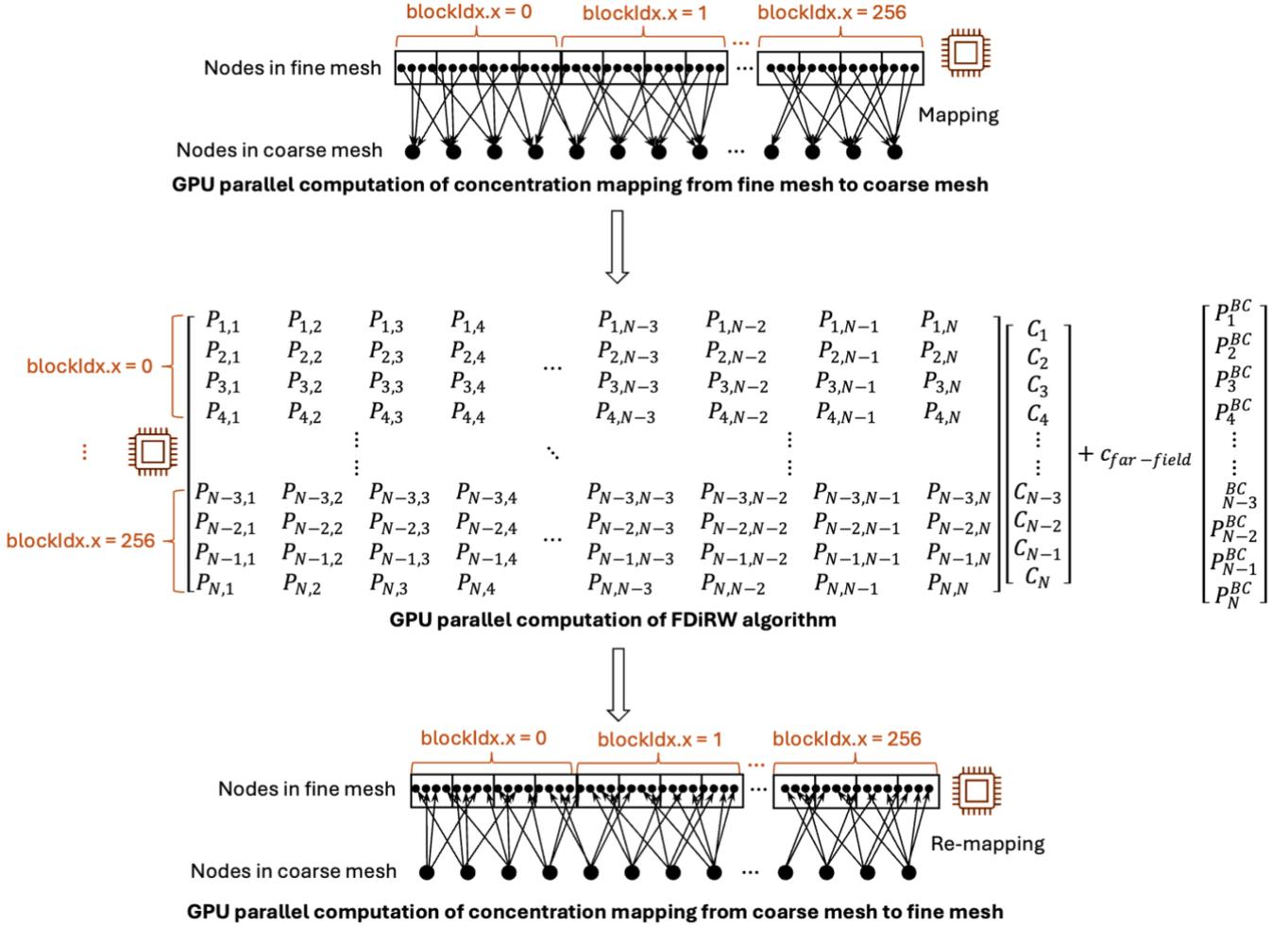

Figure 2 Illustration of GPU multi-thread parallel computation of FDiRW algorithm when solving the fast diffusion in liquid phase.

### 3.3 Mixed-precision FDiRW

The computational efficiency of FDiRW can be further enhanced by utilizing lower precision formats, such as the half-precision floating point (FP16) available since CUDA version 7.5[26]. Compared to the 64-bit double-precision FP64, which allows for approximately 15 to 17 significant digits, and the 32-bit single-precision FP32, which provides around 7 significant digits, FP16 employs only 16 bits and offers 3 to 4 significant digits. This reduced length of significant digits makes FP16 more memory-efficient and reduces the computational workload. Theoretically, FP16 can be twice as fast as FP32 and four times faster than FP64. In applications with high memory demands, the performance benefits of using lower precision can be even better.

The lower-precision computation is applicable to the FDiRW solver because the mass fractions moving across the domain, represented by the elements of the **P** matrix, do not require high precision as long as mass conservation is maintained. Since mass conservation is assured in the far-field liquid concentration calculations in Eq. (7), the low-precision computation of FDiRW superposition does not impact mass conservation. While FP16 precision might affect the exact mass fractions, 3 to 4 significant digits are typically sufficient. However, careful consideration is needed during the matrix-vector multiplication step.

When calculating concentrations at each node, the mass fractions from each node are accumulated, forming a very large sum number compared to each individual mass fraction value. This process involves adding very small numbers to much larger numbers. In such cases, the limited precision of FP16 is not adequate. Therefore, we propose a mixed-precision FDiRW approach. In this approach, accumulation operations are performed with FP32 precision, while intensive multiplication operations use FP16 precision, as illustrated in Figure 3.

Specifically, the addition-dominant 'mapping' step is executed in FP32 precision. The resulting concentration matrix, FP32 **C**, is then converted to FP16 precision for the subsequent matrix-vector multiplication step. The coefficient matrix **P** is stored in FP16 precision, so multiplication operations use FP16 precision. Then, the products of **P** and **C** are converted back to FP32 for the accumulation/addition operations. This process is also applied to the computation of boundary effect terms. The resultant concentration array **C** is obtained in FP32 precision and then converted to FP64 and assigned to the corresponding liquid nodes in the fine mesh during the 're-mapping' step. The updated FP64 concentration **c** of liquid nodes is then used in other high-precision computation steps.

This mixed-precision FDiRW algorithm aims to optimize efficiency while maintaining the desired level of accuracy.

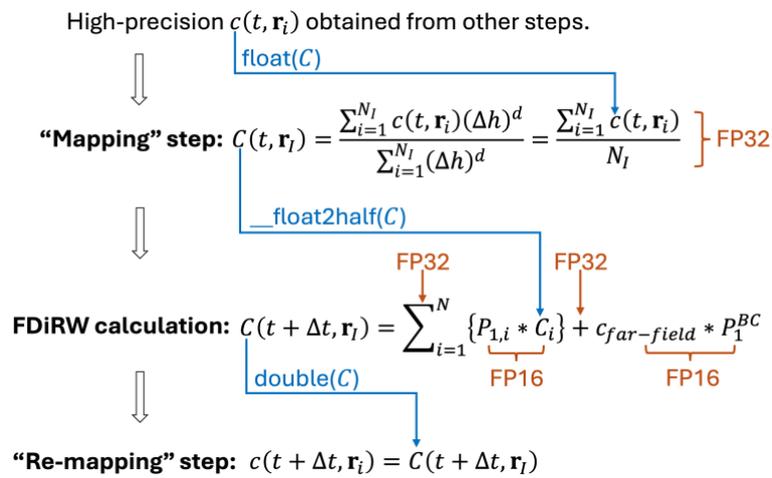

Figure 3  Illustration of mixed-precision FDiRW algorithm.

### 3.4  Flowchart

Figure 4 presents the overall flowchart of the proposed GPU-accelerated, mixed-precision integrated framework for solving the strongly inhomogeneous diffusion problem across a porous wasteform particle.

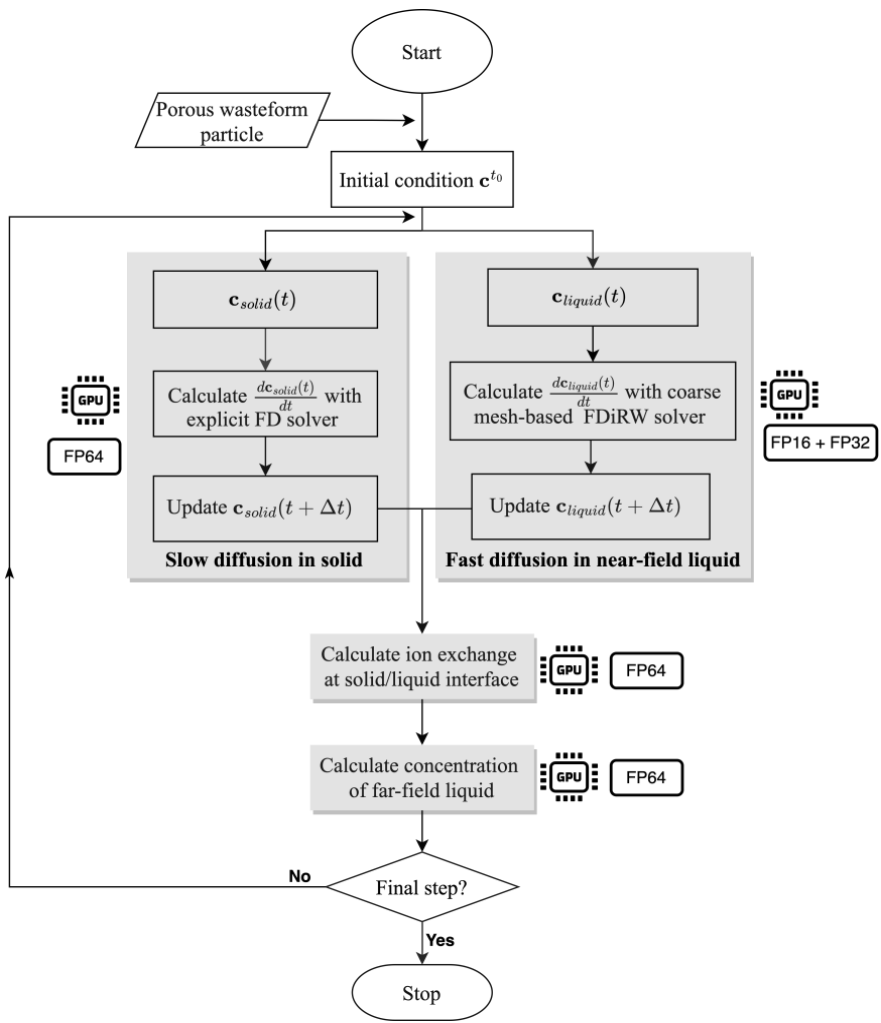

Figure 4 Flowchart of integrated solver for solving the strongly inhomogeneous diffusion problem across a porous wasteform particle.

## 4. Verification and Efficiency

### 4.1 GPU acceleration

The accuracy of the GPU-accelerated FDiRW solver is first evaluated by comparing its results to those obtained from validated CPU computations[12,20], as presented in Figure 5. The GPU solver precisely reproduces the concentration gradient in the near-field liquid, matching the corresponding CPU solutions (as shown in Figure 5a and 6b). It also yields identical concentration contours in the solid phase (Figure 5c and 6d). The effectiveness of the GPU-accelerated FDiRW solver is further assessed by examining the radionuclide absorption kinetics. We monitor the total concentration of radionuclides absorbed by the wasteform particle and those remaining in the liquid phase. The kinetic curves obtained from GPU-accelerated FDiRW solver and the benchmarking Finite Difference solver on both CPUs and GPUs are compared in Figure 6. The results confirm the desired accuracy of the GPU-accelerated FDiRW solver when compared to the fine mesh-based Finite Difference solver.

After validating the accuracy of the GPU-accelerated FDiRW solver for solving fast diffusion in the liquid phase, we examined the computational efficiency enhancement from GPU acceleration, as presented in Figure 7. In all the following comparisons, we used 1-core AMD EPYC 7763 CPU with 16GB memory storage and 1-core NVIDIA V100 GPU with 16GB memory storage and CUDA version 12.2.

As aforementioned, the simulation of the strongly inhomogeneous diffusion problem includes four main components. The efficiency performance of FDiRW solver is evaluated by comparing to the corresponding Finite Difference solver. Results in Figure 7 reveal that the computational time spent on solving fast diffusion in the near-field liquid accounts for approximately 99.8% of the total computational time when running on CPUs (Figure 7a) and GPUs (Figure 7c). This proportion reduces to 34% and 32.1% when using FDiRW solver on CPUs (Figure 7b) and GPUs (Figure 7d), respectively. Specifically, the results in Figure 7e show that, compared to the traditional FD solver, using our proposed FDiRW solver, the computational time on CPUs reduces from 909,40 seconds to 82 seconds, demonstrating an acceleration of 1109X. With GPU acceleration, the computation time for FD solver reduces from 909,40 to 757 seconds, whereas the run time for our FDiRW solver reduces from 82 seconds to 0.7 second. For both solvers, the GPU acceleration demonstrates more than 100x efficiency improvement. More impressively, from traditional FD solver on CPU to our FDiRW solver on GPU, the overall speedup is more than 100,000x.

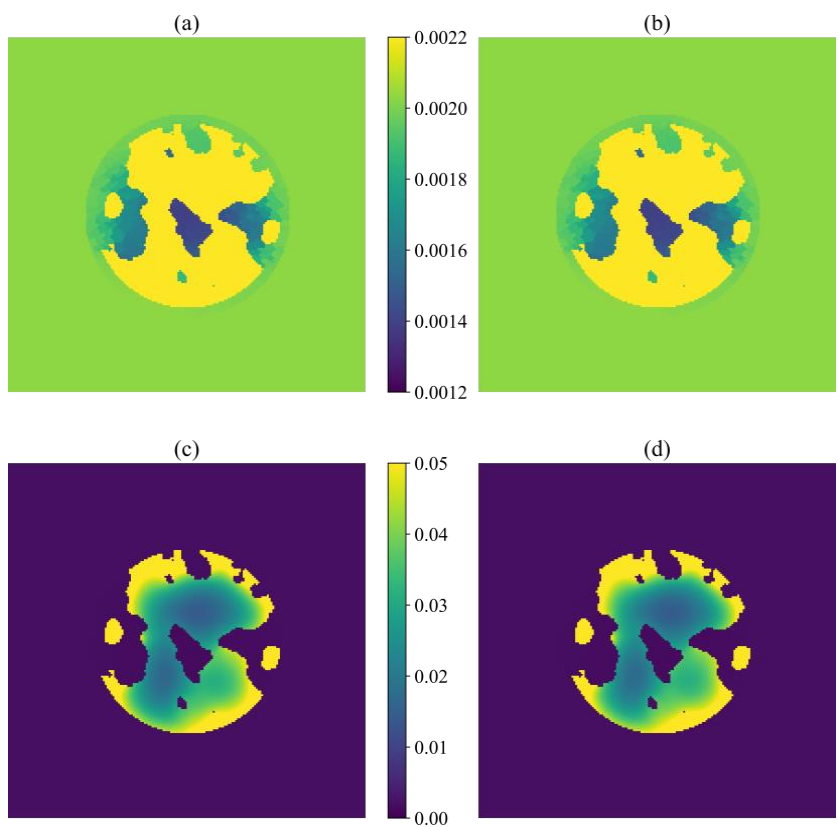

Figure 5 Screenshots of concentration fields obtained from simulations at $t = 0.5\ sec$. (a) Concentration in near-field liquid running on CPUs; (b) Concentration in near-field liquid running on GPUs; (c) Concentration in solid phase running on CPUs; (d) Concentration in solid phase running on GPUs.

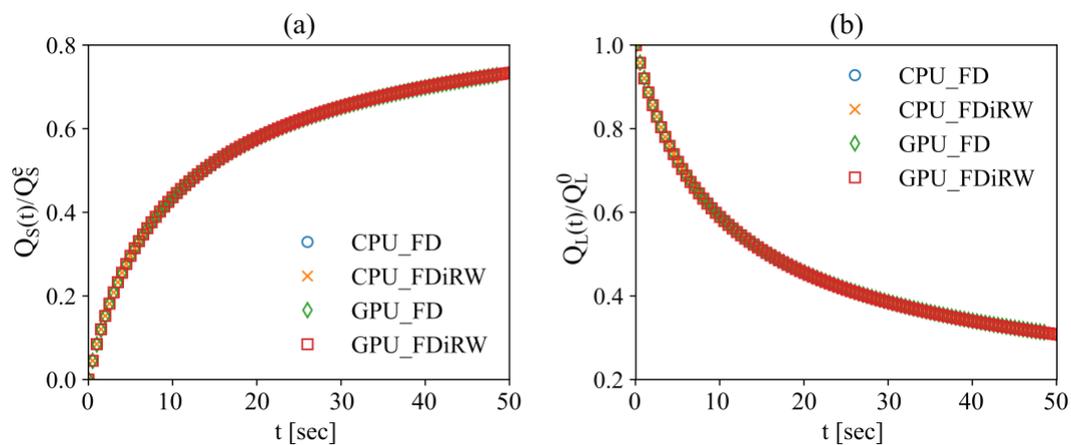

Figure 6 Radionuclide absorption kinetics obtained from simulations. $Q_S(t)$ and $Q_L(t)$ are the instant total concentration of radionuclides in solid and liquid, respectively. The absorption capacity of wasteform particles $Q_S^e = c_S^{eq} N_S$, with $N_S$ being the total number of solid nodes. The total concentration initially existed in the solution is calculated as $Q_L^0 = c_L^0 N_L$, with $N_L$ being the total number of liquid nodes in the near-field liquid phase.

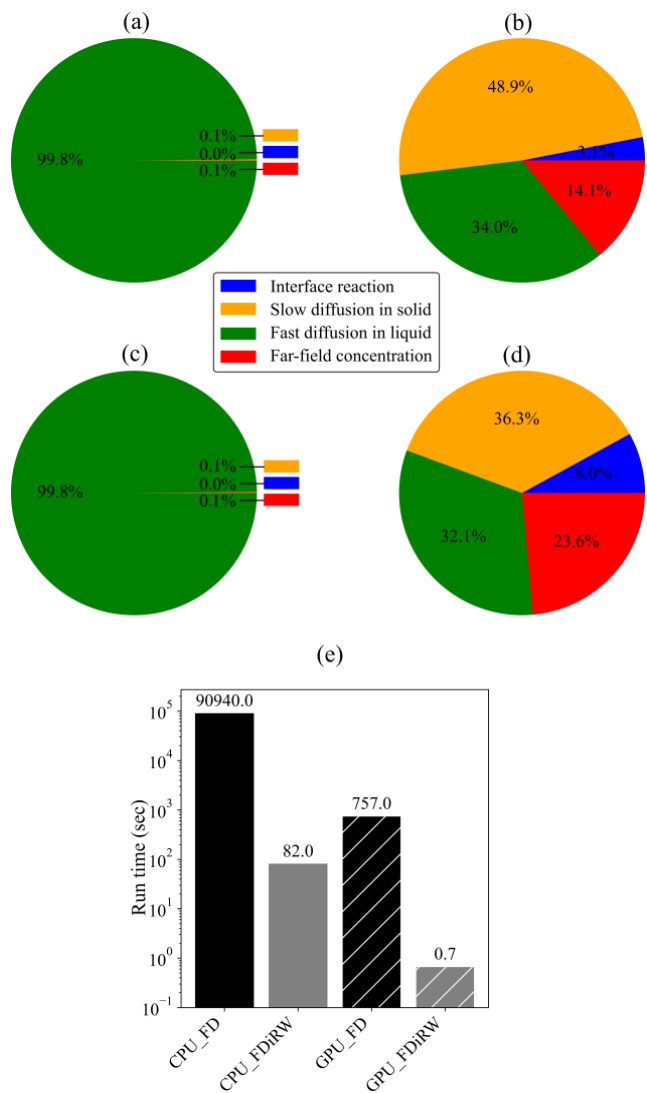

Figure 7 Comparison of run times of FD vs FDiRW solvers for solving fast diffusion in the near-field liquid within the physical time $t = 0.5 \ sec$. (a) Computation time proportions of FD running on CPUs; (b) computation time proportions of FDiRW running on CPUs; (c) Computation time proportions of FD running on GPUs; (d) computation time proportions of FDiRW running on GPUs; (e) Run time comparison of FD vs FDiRW running on CPUs and GPUs.

### 4.2 Mixed-precision FDiRW solver

More speedups can be achieved through mixed-precision computation on GPUs. However, to ensure numerical accuracy with mixed precision, thorough examinations are required. The accuracy of the FDiRW solver, using FP32, mixed FP32/FP16, and FP16, was evaluated by comparing results against the high-precision FP64 baseline. Figure 8 illustrates the concentration contours in both the liquid and solid phases obtained from FDiRW solvers with the varying precisions. A significant discrepancy is observed between the FP16 and FP64 results (Figure 8d), whereas the FP32 and mixed FP32/FP16 results closely match the FP64 results. This finding is highlighted by the error fields of the FDiRW solver using lower-precision computation, as depicted in Figure 9. The error scale of the FP16 results in the liquid phase (Figure 9 c1) is three to four orders of magnitude larger than the mixed FP32/FP16 (Figure 9 b1) and FP32 (Figure 9 a1) results. Given the concentration value scale in the liquid phase shown in Figure 8 a1, the relative error of the FP16 results reaches approximately 10%, resulting in a notably different concentration field in the solid phase (Figure 8 d2) and thus yielding inaccurate numerical solutions. Conversely, the mixed FP32/FP16 solution demonstrates a relative error of approximately 0.001%, which is close to the FP32 result of 0.0001%, and is acceptable for accuracy considerations in the domain applications.

Figure 10 further explores how errors from lower-precision computations evolve over time. Figure 10a highlights its impact on radionuclide absorption kinetics, showing that both FP32 and mixed FP32/FP16 yield results identical to FP64, while FP16 significantly affects the kinetics curve, leading to unreliable results. Figure 10b and Figure 10c displays the temporal changes of absolute errors and relative errors, respectively. FP16 results in unacceptably high errors on the order of $10^{-2}$, whereas FP32 and mixed FP32/FP16 achieve errors on the order of $10^{-6}$ and $10^{-5}$, respectively. Additionally, the numerical errors caused by lower-precision computation do not continue to grow over iterations in the computation. These findings demonstrate desired accuracy of the mixed FP32/FP16 computation throughout the process and suggest the inadequacy of FP16.

The inadequacy of FP16 in the FDiRW solver is attributed to the large range of element magnitudes in the coefficient matrices $\mathbf{P}$ and $\mathbf{P}_{BC}$, as shown in Figure 11. Adding values of significantly different scales—such as a very small number to a very large number—requires sufficient significant digits. FP16 precision, with only 3-4 significant figures, is insufficient, leading to large numerical errors. Furthermore, even when elements in the coefficient matrix are of similar scales, matrix-vector multiplication involves summing the products of corresponding elements ($P_{ij}$ and $C_j$), which accumulates to a scale several orders of magnitude larger than the individual addends. This summation demands long significant digits for desired accuracy. Therefore, using FP32, which provides longer significant digits for addition operations, substantially improves accuracy compared to FP16.

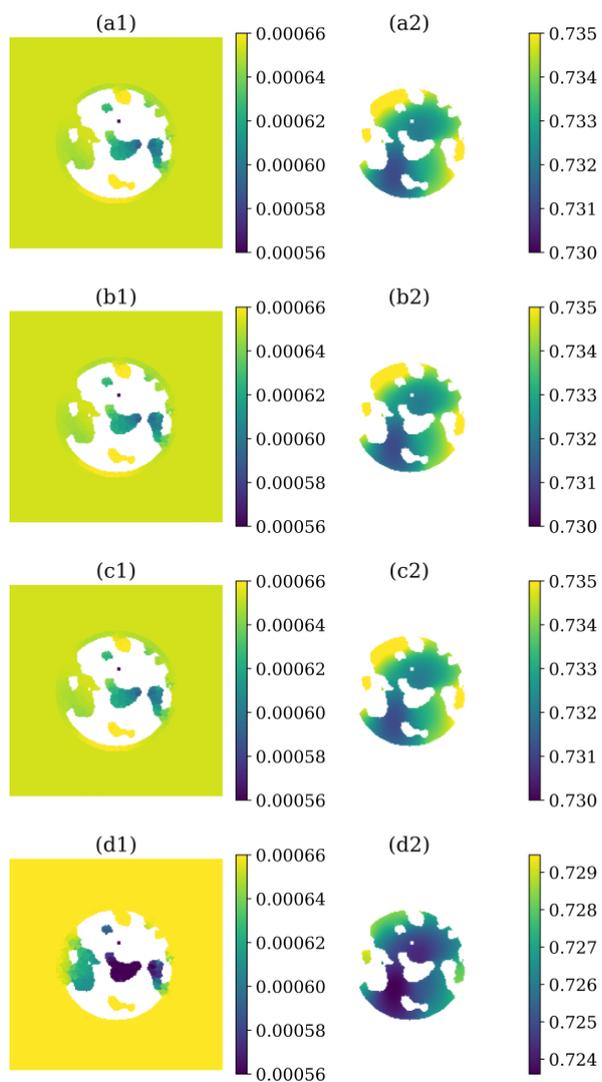

Figure 8 Concentration contours in liquid phase (Left) and solid phase (Right) obtained from FP64 FDiRW solver at physical time t = 50 sec. (a) FP64, (b) FP32, (c) FP32/FP16, (d) FP16.

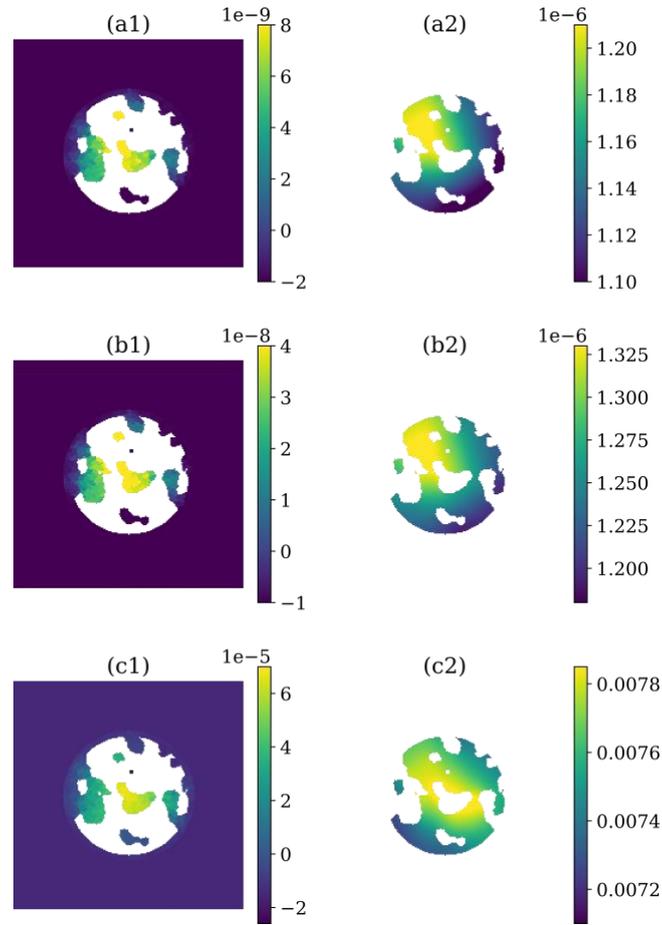

Figure 9 Contours of concentration discrepancy of low-precision solution from FP64 solution in liquid phase (Left) and solid phase (Right) at physical time t = 50 sec. (a) FP32, (b) FP32/FP16 (c) FP16. The concentration discrepancy is calculated as the absolute error of concentration obtained from lower-precision computation relative to the FP64 results.

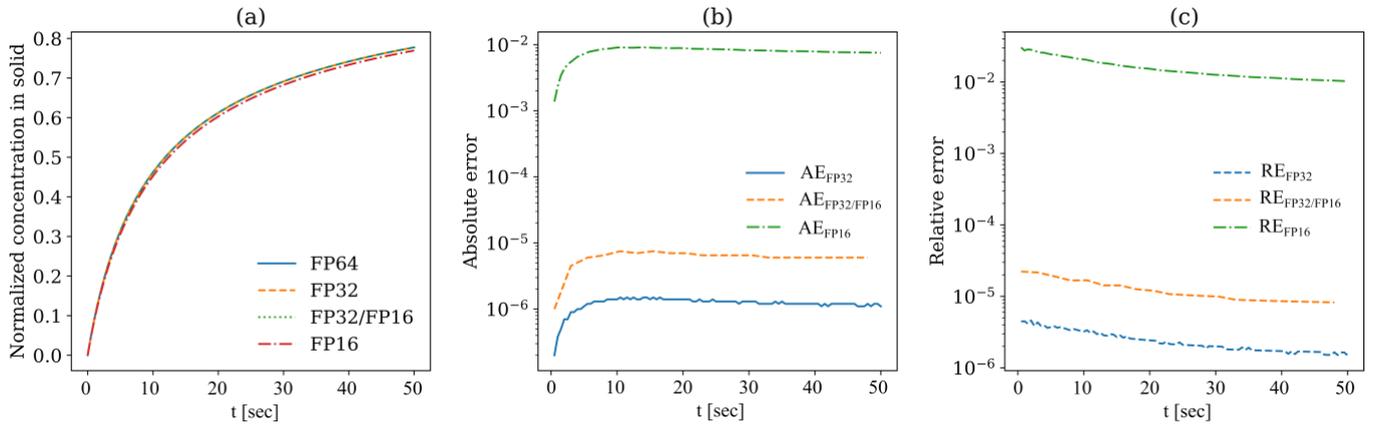

Figure 10 Accuracy comparison of different precisions based FDiRW solver. The normalized concentration $\bar{c}_S(t)$ in solid presented in (a) is the solid nodes' averaged concentration normalized by $c_S^{eq}$. The absolute error in (b) is calculated as $AE_{FP*} = |\bar{c}_{S_{FP64}} - \bar{c}_{S_{FP*}}|$, and the relative error in (c) is calculated as $RE_{FP*} = |\bar{c}_{S_{FP64}} - \bar{c}_{S_{FP*}}|/\bar{c}_{S_{FP64}}$, where $\bar{c}_{S_{FP64}}$ and $\bar{c}_{S_{FP*}}$ are the normalized concentration in solid computed with FP64 precision and other precisions, respectively.

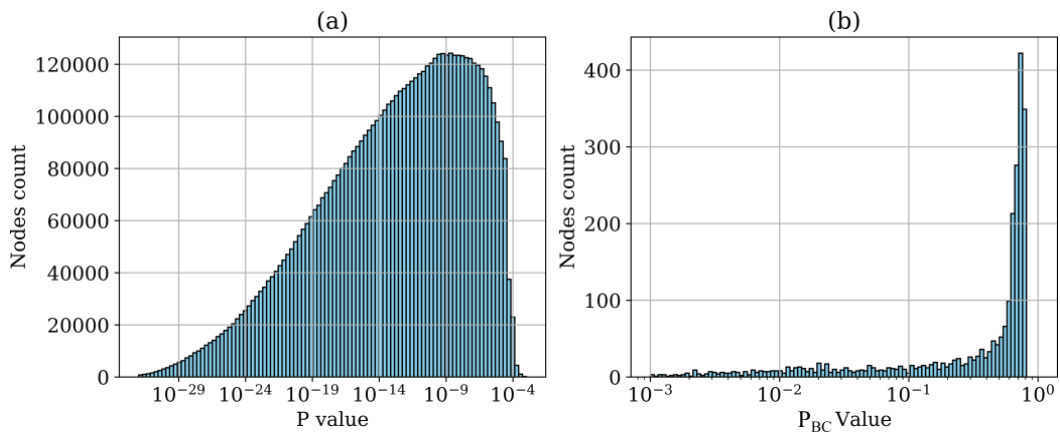

Figure 11 Distribution of element values of coefficient matrix $\mathbf{P}$ (a) and $\mathbf{P}_{BC}$ (b).

After verifying the accuracy of the proposed mixed-precision FDiRW solver, we compared its computational efficiency to other precision configurations as listed in Table 2 and plotted in Figure 12. This comparison comprises six cases: FP64, FP32, mixed FP32/FP16, FP16, and two additional cases using the cuBLAS library for matrix-vector multiplications under FP64 and FP32 precisions. The numerical tests reveal that, compared to the baseline FP64, FP32 achieves a 1.55x speedup, while FP16 achieves a 1.91x speedup. These results fall short of the theoretical 2x and 4x speedups due to the overhead associated with precision conversions when integrated with alternative high-precision computation in other components. The mixed FP32/FP16 configuration achieves a 1.75x speedup with respect to FP64, offering a nice balance between accuracy and efficiency that makes it the best choice overall. Although faster, FP16 cannot maintain desired accuracy.

Notably, using cuBLAS library packages with optimizations for computation and memory storage results in 1.48x and 2.40x speedups for FP64 and FP32, respectively. Among all six cases, FP32 with cuBLAS provides the best efficiency with desirable accuracy. This finding underscores the importance of optimizing GPU computation and memory storage for the proposed mixed FP32/FP16 FDiRW solver. While such optimizations are beyond the scope of this study, they suggest significant potential for further enhancing efficiency.

Table 2 GPU run time of different precisions of FDiRW solvers.

| Precisions | GPU run time [sec] | Normalized GPU run time |
|---|---|---|
| FP64 | 75.3 | 1.0 |
| FP32 | 48.7 | 1/1.55 |
| FP32/FP16 | 43.1 | 1/1.75 |
| FP16 | 39.5 | 1/1.91 |
| FP64_cublasDgem | 50.8 | 1/1.48 |
| FP32_cublasSgem | 31.4 | 1/2.40 |

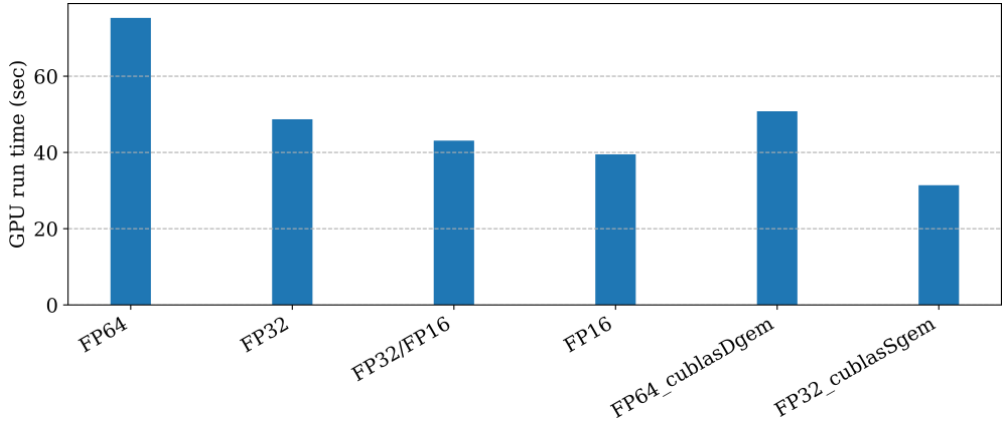

Figure 12 Efficiency comparison of different precisions of FDiRW solvers.

### 4.3 Applications to larger particles

The GPU-accelerated mixed-precision FDiRW solver is then applied to varying sizes of porous particles, as shown in Figure 13, to investigate its scalability to larger scale systems as the dimension of the coefficient matrix increases. Larger particles exhibit a more prominent concentration gradient inside the solid particles due to that large particles contain regions with larger solid volumes, and thus it takes longer for radionuclide ions to reach those areas and achieve absorption capacity, resulting in slower kinetics compared to smaller particles. This trend is illustrated in Figure 14, where the normalized concentration in large particles converges to 1.0 at a slower rate than in small particles.

The scalability of the proposed FDiRW solver, when applied to larger particles or multiple-particle systems, is presented in Figure 15 with further details provided in Table 3. It shows that the computation complexity of FDiRW is O($N_L$), as illustrated in Figure 15a, keeps consistent across all precision conditions. To analyze the computational complexity, we examine the floating-point operations (FLOPs) of the FDiRW solver breaking down into three steps: the 'mapping' step, FDiRW calculation step, and 'remapping' step, as given in Figure 3. In this FLOPs analysis, we focus on the computationally intensive multiplication operations while ignoring the addition operations considering their relatively lower computational burden. During the 'mapping' step, each liquid node's concentration is mapped to its corresponding representative node, resulting in $N_L$ FLOPs. Similarly, in the re-mapping step, the concentration of representative nodes is remapped back to each corresponding liquid nodes, also requiring $N_L$ operations. In the FDiRW calculation step, the FLOPs required for matrix-vector multiplication in Eq. (10) is $N \times N + N$, including $N \times N$ FLOPs for computing $\mathbf{P}|_{N \times N} \mathbf{C}|_{N \times 1}$, and $N$ FLOPs for computing $P_{BC}|_{N \times 1} c_{far-field}$. Hence, the total FLOPs for the three-step FDiRW solver is $N(N + 1) + 2 * N_L$. Depending on the relative values of the total number of liquid nodes $N_L$ and the total number of representative nodes $N$, the final computational complexity can be $O(N^2)$ if dominated by the FDiRW calculation step, or $O(N_L)$ if dominated by the 'mapping' and 'remapping' steps. The numerical tests results shown in Figure 15 shows that the computation is essentially dominated by the 'mapping' and 'remapping' steps, yielding a computation complexity of $O(N_L)$. This may be attributed to the frequent access to the index memory of representative nodes corresponding to each node in the near-field liquid during the 'mapping' and 'remapping' steps, causing a significant memory overhead relative to computation.

Figure 15c compares the efficiency gain of mixed-precision computation when applied to different model sizes. The results indicate that the FDiRW solver consistently benefits from mixed-precision computation. A closer examination shows that the efficiency gain slightly increases as the model size grows. This trend can be explained by the fact that larger models require more substantial computational loads, thereby reducing the relative impact of overheads from other sources. Consequently, the efficiency enhancement achieved through mixed-precision becomes more pronounced and noticeable with scaled model sizes.

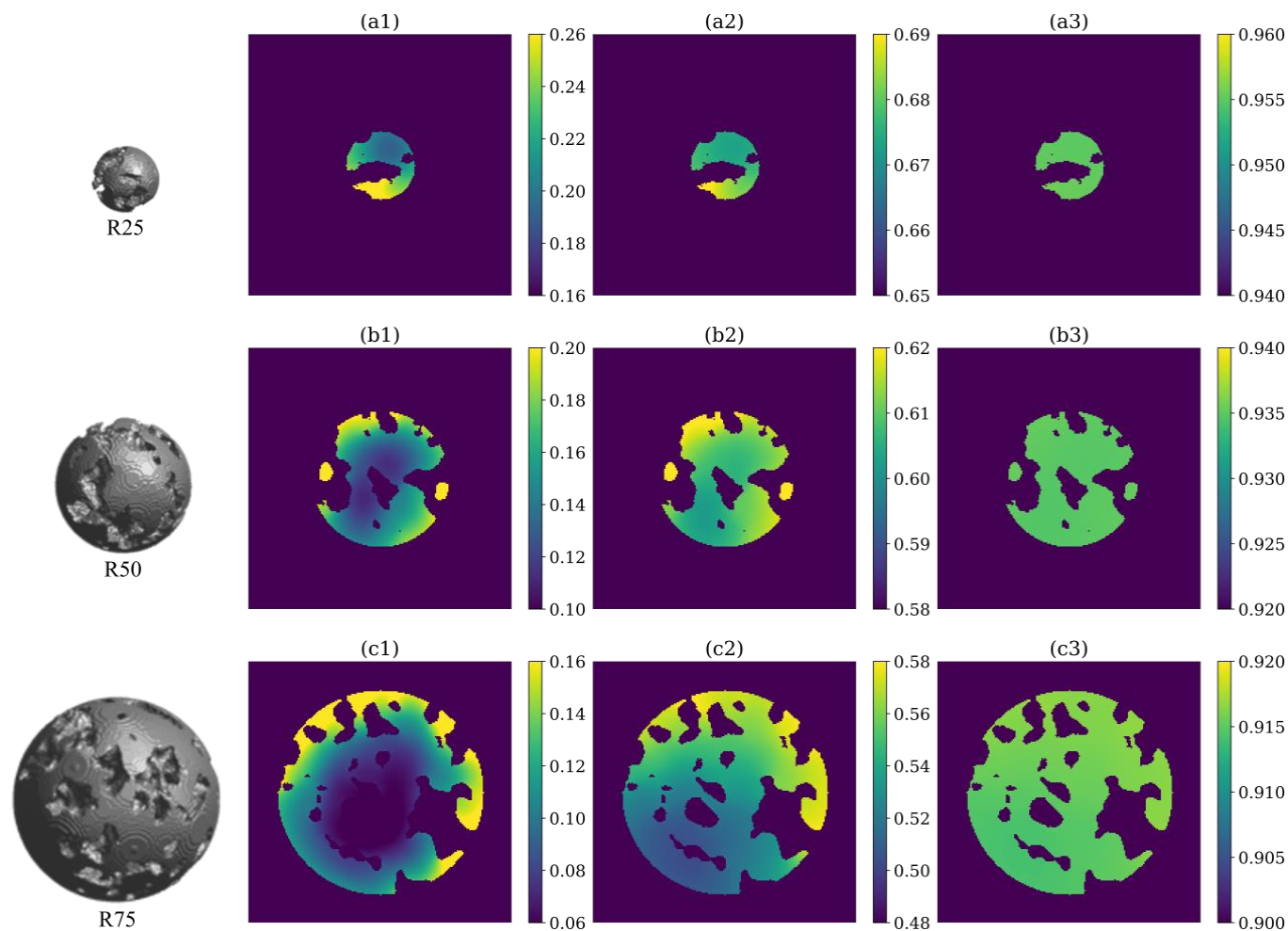

Figure 13  Geometries of three different sizes of porous wasteform particles obtained from phase field modeling. The three particles denoted as R25, R50, R75 have a radius of $25\Delta h, 50\Delta t, 75\Delta h$, respectively. Note that, the colorbar range is chosen to exhibit the concentration gradient in solid phase, and the dark color in liquid phase does not mean the liquid has a value of the lowest value of the colorbar range. The concentration in liquid phase is very small and approaches to zero over time.

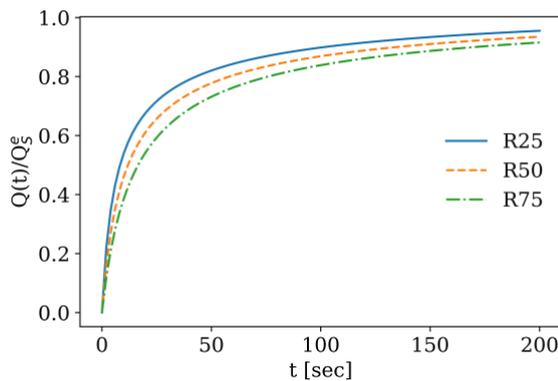

Figure 14 Radionuclide absorption kinetics of different sizes of particles obtained from FD-FDiRW.

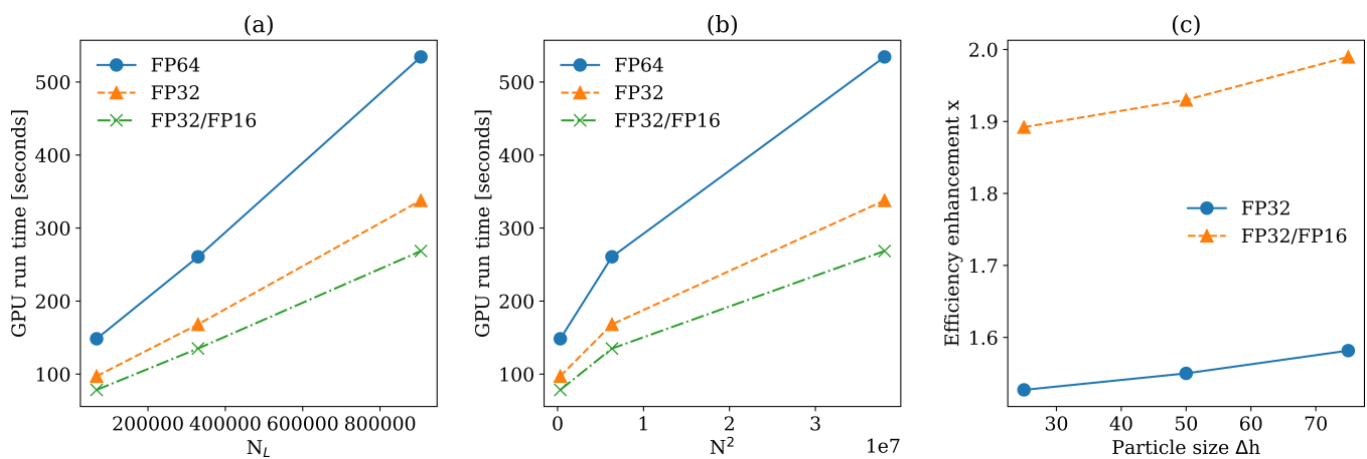

Figure 15 Scalability of FDiRW solver when applied to large porous wasteform particles. (a) GPU run time *vs* $N_L$; (b) GPU run time *vs* $N^2$. (c) Efficiency gain by using lower precisions. GPU run time is linearly proportional to $N_L$, meaning the computational complexity of FDiRW is O($N_L$). In (c), the efficiency enhancement is calculated by as the ratio of computational time of FP64 to that of lower precisions.

Table 3 Total number of nodes and representative nodes in near-field liquid of the three particles R25, R50, R75, and the corresponding GPU run time in FDiRW computation.

|  |  | **R25 particle** | **R50 particle** | **R75 particle** |
|---|---|---|---|---|
| $N_L$ |  | 67,207 | 329,404 | 905,426 |
| $N_{rl}$ |  | 561 | 2,515 | 6,168 |
| $N_S$ |  | 45,724 | 366,713 | 1,238,185 |
| $V_L^{far-field}$ |  | $2.27 \times 10^{-11}$ [mL] | $1.80 \times 10^{-10}$ [mL] | $6.16 \times 10^{-10}$ [mL] |
| **GPU run time [sec]** | FP64 | 148 | 260 | 534 |
|  | FP32 | 97 | 168 | 338 |
|  | FP32-FP16 | 78 | 135 | 268 |

## 5. Conclusions and Discussions

In conclusion, this study introduces and validates a GPU-accelerated mixed-precision configuration for the FDiRW solver, tailored to address computationally intensive multiphase diffusion problems with strongly inhomogeneous diffusivity. Through a single GPU's parallel acceleration, the solver achieves 117x speedup compared to a single-core CPU for a medium-sized model of 192×192×192. Furthermore, the proposed mixed FP32/FP16 configuration provides an additional 1.75x speedup in GPU computation, achieving an optimal balance between accuracy and efficiency. In terms of scalability, numerical tests indicate that the computational complexity of the solver is O($N_L$), where $N_L$ is the total number of nodes in computation. This demonstrates that the FDiRW approach is effective and cost-efficient for tackling strongly inhomogeneous problems in scientific and engineering applications.

Despite its significant advantages in computational efficiency, the proposed FDiRW approach has its own limitations. A notable one stems from the assumption of a fixed diffusion domain, as the FDiRW solver relies on a preconditioning coefficient matrix **P**. In other words, FDiRW is not well-suited for scenarios where the porous structure evolves over time, leading to changes of domain where the fast diffusion takes place. In such more dynamic environments, the coefficient matrix **P** needs to be updated continuously to reflect the changing configuration of the fluid domain. These frequent updates to the **P** matrix can significantly undermine the potential efficiency gains from the FDiRW method. Furthermore, these frequent updates of **P** also complicates the implementation process.

Nevertheless, the FDiRW solver has a wide range of applications involving strongly inhomogeneous diffusion processes with fixed or slowly evolving microstructures. The combination of the efficient FDiRW algorithm and the substantial acceleration provided by GPU mixed-precision computation has significantly reduced the computational time to just a few minutes on a single GPU. This improvement paves the way for the investigation of more realistically large systems, including those with multiple particles or larger particles. For instance, in experimental settings, particle sizes are typically in the range of tens to hundreds of microns, which is 2 to 3 orders of magnitude larger than the particle sizes studied in this paper. The GPU-accelerated FDiRW solver enables the study of these larger, more complex real-world systems, making it feasible to explore and analyze diffusion processes in environments that were previously too computationally intensive to handle.

Beyond solving the fast diffusion equation, the FDiRW solver can be extended to include dynamical flow which accelerates the transport of mobile species in liquid phase. While our focus in this study is on solving the diffusion equation, the FDiRW approach is capable of solving the diffusion-convection equation when calculating the elements of the coefficient matrix **P**. This extension is feasible as long as the dynamical flow can be regarded as a steady process and independent of other processes. This means that the flow's influence on diffusion can be integrated into the coefficient matrix **P** without requiring continuous updates during the computation. Once the coefficient matrix **P** is determined, incorporating dynamical flow does not alter the matrix-vector multiplication formulation in Eq. (8) nor the mapping/remapping steps. Therefore, adding dynamical flow is straightforward and does not increase the computational demands of the integrated numerical framework. This allows for a seamless extension of the FDiRW solver to more complex scenarios involving both diffusion and convection, enhancing its versatility and applicability in various scientific and engineering contexts.

**CRediT authorship contribution statement**

**Zirui Mao**: Investigation, Methodology, Visualization, Writing - original draft. **Shenyang Hu**: Writing - review & editing. **Ang Li**: Funding acquisition, Discussion, Writing - review & editing.

**Acknowledgements**

The work described in this article was performed by Pacific Northwest National Laboratory, which is operated by Battelle for the U.S. Department of Energy under Contract DE-AC05-76RL01830. This material is based upon work supported by the U.S. Department of Energy, Office of Science, Office of Advanced Scientific Computing Research, ComPort: Rigorous Testing Methods to Safeguard Software Porting, under Award Number 78284. This work was also partially supported by the Center for Hierarchical Waste Form Materials, an Energy Frontier Research Center funded by the U.S. Department of Energy, Office of Science, Basic Energy Sciences under Award No. DE-SC0016574. This research used resources of the National Energy Research Scientific Computing Center (NERSC), a U.S. Department of Energy Office of Science User Facility located at Lawrence Berkeley National Laboratory, operated under Contract No. DE-AC02-05CH11231.